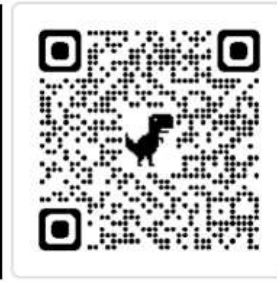

# Versatile Stochastic Two-Sided Platform Models

SONG-KYOO (AMANG) KIM

**ABSTRACT**

This paper deals with the alternative mathematical modeling of the two-side platform. Two-sided platforms are specific multi-sided platforms that bring together two distinct groups of a model. The stochastic modeling by adapting various innovative mathematical methods including the first exceed theory and the stochastic pseudo-game theory has been applied for describing a two-sided platform more properly. A stochastic pseudo-game model is newly introduced to solve the two-sided platform more effectively. Analytically tractable results for operation thresholds for maximizing profits are provided and it also delivers the optimal balance of a two-sided platform. The paper includes how these innovative models are applied into various two-sided market situations. Additionally, users could conduct these multi-sided models to real business developments and the case practices of these unique models shall help the readers who want to find recommendations of their business situations easily even without having any mathematical background.

## I. INTRODUCTION

Multi-sided platforms bring together two or more distinct but interdependent groups of customers, normally described as business-to-business (B2B) and business-to-customer (B2C). Two-sided platforms have proliferated rapidly with the Internet and this has led to the development of new business models to monetize innovative value propositions not only in online but also offline markets [1-4]. Multi-sided platforms are an important business phenomenon that has proliferated with the rise of information technology and the Internet. Two-sided platforms are specific multi-sided platforms that bring together two distinct but interdependent groups of customers [5-13]. A two-sided market involves two different user groups whose interactions are usually enabled over an open platform. The two sides of the market represent the two primary sets of economic agents while the platform acts as an enabler or catalyst for bringing together the two distinct set of economic agents [14-15]. In such markets, both categories of economic agents have to be present on the platform in right proportions to create enough value to both sides and thereby accelerate and/or sustain the platform. Hence, in such market systems, one side's participation depends on the value created by presence of the other side over the platform. This is termed as the cross-side network effect or network externalities [16-18].

Although there are several researches to use mathematical models to describe multi-sided platform, the stochastic modeling by adapting various innovative mathematical methods including the first exceed theory and the pseudo-stochastic game theory has been applied for describing a two-sided platform more properly.

## II. Stochastic Modelling For Two-sided Platform

The innovative stochastic model has been adapted into the two-sided platform to analyze the behaviors of both sides. This analytical model consists of two sides (i.e., the customer and the supplier sides). This explicit function (Theorem 1) gives the predicted moment of one step prior to hit the platform capacity by customers (simply called player A). The game theory framework makes complicated game types of problems simple to understand situations more clearly. Therefore, stochastic models based on a game framework have been widely studied [19, 20] and these models are applied for various areas including the Blockchain [21-23] and the business strategy [24-26]. The recent game based stochastic model is adapted into the stock market exchanges [20]. A stochastic pseudo-game is not a stochastic game model but only adapts the stochastic game framework. The stochastic model in this research is a stochastic pseudo-game because this model does not describe a game type of situation but adapts the stochastic game framework for designing the two-sided platform.

### 2.1. Stochastic Pseudo-Game

The antagonistic two-person pseudo-game (called "A" and "B") describes the behavior of a two-sided platform between a customer (player A) and supplier (player B) sides. Both input flows to fill the capacity of the platform from either from a customer or a supplier sides. Let $(\Omega, \mathcal{F}(\Omega), P)$ be probability space $\mathcal{F}_a, \mathcal{F}_b, \mathcal{F}_\tau \subseteq \mathcal{F}(\Omega)$ be independent $\sigma$-subalgebras. Suppose:

$$\boldsymbol{A} := \sum_{k \geq 0} X_k \varepsilon_{s_k}, \; s_0(=0) < s_1 < s_2 < \cdots, \text{ a.s.} \tag{2.1}$$

$$\boldsymbol{B} := \sum_{j \geq 0} Y_j \varepsilon_{t_j}, \; t_0(=0) < t_1 < t_2 < \cdots, \text{ a.s.} \tag{2.2}$$

are $\mathcal{F}_a$- and $\mathcal{F}_b$-measurable marked Poisson processes ($\varepsilon_w$ is a point mass at $w$) with respective intensities $\lambda_a$ and $\lambda_b$ which are related with the input behaviors of both sides in the two-sided platform. These processes represent the customer flow of player A (a customer side) and the number of suppliers of player B (a supplier side). Player A checks the number of customer input at times $s_1, s_2, \ldots$ and sustain respective customers $X_1, X_2, \ldots$ formalized by the process $\boldsymbol{A}$. The number of suppliers to player B are described by the process $\boldsymbol{B}$ similarly. Player B counts the number of suppliers (or contents) to cover the customers in the platform. The processes $\boldsymbol{A}$ and $\boldsymbol{B}$ are specified by their transforms

$$\mathbb{E}\left[z^{\boldsymbol{A}(s)}\right] = e^{\lambda_a s(z-1)}, \; \mathbb{E}\left[g^{\boldsymbol{B}(t)}\right] = e^{\lambda_b t(g-1)}. \tag{2.3}$$

The two-sided platform is monitored at random times in accordance with the point process in the two-side platform:

$$\boldsymbol{T} := \sum_{i \geq 0} \varepsilon_{\tau_i}, \ \tau_0(\,>0)), \tau_1, \ldots, \qquad (2.4)$$

which is assumed to be delayed renewal process. If

$$(A(t), B(t)) := \boldsymbol{A} \otimes \boldsymbol{B}([0, \tau_k]), \ k = 0, 1, \ldots, \qquad (2.5)$$

forms an observation process upon $\boldsymbol{A} \otimes \boldsymbol{B}$ embedded over $\boldsymbol{T}$, with respective increments

$$(X_k, Y_k) := \boldsymbol{A} \otimes \boldsymbol{B}([\tau_{k-1}, \tau_k]), \ k = 1, 2, \ldots, \qquad (2.6)$$

and

$$X_0 = A_0, \ Y_0 = B_0. \qquad (2.7)$$

The observation process could be formalized as

$$\boldsymbol{A}_\tau \otimes \boldsymbol{B}_\tau := \sum_{k \geq 0} (X_k, Y_k) \varepsilon_{\tau_k}, \qquad (2.8)$$

where

$$\boldsymbol{A}_\tau = \sum_{i \geq 0} X_i \varepsilon_{\tau_i}, \ \boldsymbol{B}_\tau = \sum_{i \geq 0} Y_i \varepsilon_{\tau_i}, \qquad (2.9)$$

and it is with position dependent marking and with $X_k$ and $Y_k$ being dependent with the notation

$$\Delta_k := \tau_k - \tau_{k-1}, \ k = 0, 1, \ldots, \ \tau_{-1} = 0, \qquad (2.10)$$

and

$$\gamma(z, g) = \mathbb{E}\left[z^{X_k} \cdot g^{Y_k}\right], \ g > 0, \ z > 0. \qquad (2.11)$$

By using the double expectation,

$$\gamma(z, g) = \delta(\lambda_a(1 - z) + \lambda_b(1 - g)), \qquad (2.12)$$

and

$$\gamma_0(z, g) = \mathbb{E}\left[z^{A_0} g^{H_0}\right] = \delta_0(\lambda_a^0(1 - z) + \lambda_b^0(1 - g)), \qquad (2.13)$$

where

$$\delta(\theta) = \mathbb{E}\left[e^{-\theta \Delta_1}\right], \ \delta_0(\theta) = \mathbb{E}\left[e^{-\theta \tau_0}\right], \qquad (2.14)$$

are the magical transform of increments $\Delta_1, \Delta_2, \ldots$. The stochastic process $\mathcal{A}_\tau \otimes \mathcal{B}_\tau$ describes the status of the two-sided platform which has a connection between both sides and the supplier side is usually a dominant side. The attraction factor $\alpha$ is newly proposed which indicates the power of the dominant side to seduce the customers in the submissive side. Since we consider a supplier side is dominant which impacts on the input of customers, the customer input rate $\lambda_a$ becomes the function of $b$ and it could be assigned as follows:

$$\lambda_a(b) = \lambda_a \alpha \cdot \left(\frac{b}{\widetilde{\delta}}\right), \alpha \geq 1, b \leq M, \tag{2.15}$$

where $\widetilde{\delta} = \mathbb{E}[\Delta_1]$ and $\alpha$ is the ratio between a supplier and a customer sides and this ratio is called the attraction factor which depends on the status of the dominant side. The stochastic process $\boldsymbol{A}_\tau \otimes \boldsymbol{B}_\tau$ is completed when on the $k$-th observation epoch $\tau_k$, the collateral customers to player A exceeds the capacity of the platform $M$. To further formalize the game, the *exit index* is introduced:

$$\mu := \inf\left\{k : A_k = \sum_{k \geq 0} X_k \geq M, X_0 = A_0\right\}, \tag{2.16}$$

$$\nu := \inf\left\{j : B_j = \sum_{j \geq 0} Y_j \geq M, Y_0 = B_0\right\}. \tag{2.17}$$

Since, the platform shall be terminated at the time $\tau_\mu$, the number of customers may reach the capacity of the platform $M$. We shall be targeting the confined game in the view point of player A. The first passage time $\tau_\mu$ is the associated exit time from the confined game and the formula (2.6) will be modified as

$$\overline{\boldsymbol{A}_\tau} \otimes \overline{\boldsymbol{B}_\tau} := \sum_{k \geq 0}^{\mu} (X_k, Y_k) \varepsilon_{\tau_k} \tag{2.18}$$

which the path of the game from $\mathcal{F}(\Omega) \cap \{\nu < \mu\}$, which gives an exact definition of the model observed until $\tau_\mu$. The joint functional of the two-sided platform is as follows:

$$\Phi_M = \Phi_M(\xi, z_0, z_1, g_0, g_1) \tag{2.19}$$

$$= \mathbb{E}\left[\xi^\mu \cdot z_0^{A_{\mu-1}} \cdot z_1^{A_\mu} \cdot g_0^{B_{\nu-1}} \cdot g_1^{B_\nu} \mathbf{1}_{\{\mu < \nu\}}\right],$$

where $M$ indicates the platform capacity. This functional shall represent the status of the two-sided platform upon the exit time $\tau_\mu$. The latter is of particular interest, we are interested in not only the prediction of the moment of full but also one observation prior to this. The Theorem 1 establishes an explicit formula $\Phi_M$ from (2.11)-(2.13). The $\mathcal{D}$-operator and its inverse operator $\mathfrak{D}$ from the first exceed model [27, 28] have been adapted as follows:

$$\mathcal{D}_{(x,y)}[f(x,y)](u,v) := (1-u)(1-v)\sum_{x \geq 0}\sum_{y \geq 0} f(x,y)u^x v^y, \tag{2.20}$$

then

$$f(x,y) = \mathfrak{D}_{(u,v)}^{(x,y)}[\mathcal{D}_{(x,y)}\{f(x,y)\}], \tag{2.21}$$

where $\{f(x,y)\}$ is a sequence, with the inverse

$$\mathfrak{D}_{(u,v)}^{(m,n)}(\bullet) = \begin{cases} \left(\frac{1}{m!\cdot n!}\right) \lim_{(u,v)\to 0} \frac{\partial^m \partial^n}{\partial u^m \partial v^n} \frac{1}{(1-u)(1-v)}(\bullet), & m \geq 0, n \geq 0, \\ 0, & otherwise. \end{cases} \quad (2.22)$$

**Theorem 1:** the functional $\Phi_M$ of the process of (2.19) satisfies following expression:

$$\Phi_M = \mathfrak{D}_{(u,v)}^{(M,M)} \left[ \phi_0^1 - \phi_0 + \frac{\xi \cdot \gamma_0}{1 - \xi\gamma} (\phi^1 - \phi) \right]. \quad (2.23)$$

***Proof:*** we find the explicit formula of the joint function $\Phi_m$. The joint functional (2.19) of player A is as following:

$$\Phi_m(\xi, z_0, z_1, g_0, g_1) = \sum_{k=0}^{\infty} \xi^k \mathbb{E} \left[ \mathbf{1}_{\{\mu(m)=k\}} z_0^{A_{\mu-1}} \cdot z_1^{A_\mu} \cdot g_0^{B_{\nu-1}} \cdot g_1^{B_\nu} \mathbf{1}_{\{\mu<\nu\}} \right]$$
$$= \sum_{k=0}^{\infty} \sum_{j=k+1}^{\infty} \xi^k \mathbb{E} \left[ \mathbf{1}_{\{\mu(m)=k, \nu(n)=j\}} z_0^{A_{\mu-1}} \cdot z_1^{A_\mu} \cdot g_0^{B_{\nu-1}} \cdot g_1^{B_\nu} \right] \quad (2.24)$$

and, applying the operator $\mathcal{D}$ to random family $\{\mathbf{1}_{\{\mu(m)=k,\nu(n)=j\}} : m \geq 0\}$, we arrive at

$$\mathcal{D}_{(x,y)} \left[ \mathbf{1}_{\{\mu(m)=k,\nu(n)=j\}} \right](u,v) = \left( u^{A_{k-1}} - u^{A_k} \right) \left( v^{B_{j-1}} - v^{B_j} \right) \quad (2.25)$$

and, from the previous research [19, 27, 28],

$$\Psi(u,v) = \mathcal{D}_{(x,y)}[\Phi_m](u,v) = \sum_{k=0}^{\infty} D_{1k} D_{2k} \sum_{j>k} D_{3kj} D_{4kj}, \quad (2.26)$$

where

$$D_{1k} = \xi^k \mathbb{E}\left[ z_0^{A_{k-1}} \cdot z_1^{A_{k-1}} \cdot g_0^{B_{k-1}} \cdot g_1^{B_{k-1}} u^{A_{k-1}} \right]$$
$$= \begin{cases} 1, & k = 0, \\ \xi\gamma_0 (\xi \cdot \gamma)^{k-1}, & k > 0, \end{cases} \quad (2.28)$$

$$D_{2k} = \mathbb{E}\left[ g_1^{X_k} z_1^{Y_k} u^{X_k} (1 - u^{X_k}) \right]$$
$$= \begin{cases} \phi_0^1 - \phi_0, & k = 0, \\ \phi^1 - \phi, & k > 0, \end{cases} \quad (2.29)$$

$$D_{3kj} = \mathbb{E}\left[ v^{(Y_{k+1}+\cdots+Y_{j-1})} \right] = \gamma(1,v)^{j-(k+1)}, \quad (2.30)$$

$$D_{4kj} = \mathbb{E}[1 - v^{Y_j}] = 1 - \gamma(1,v), \quad (2.31)$$

and

$$\gamma := \gamma(z_0 z_1 u, g_0 g_1 v), \quad (2.32)$$
$$\gamma_0 := \gamma_0(z_0 z_1 u, g_0 g_1 v), \quad (2.33)$$

$$\phi := \gamma(z_1 u, g_1 v), \quad (2.34)$$
$$\phi_0 := \gamma_0(z_1 u, g_1 v), \quad (2.35)$$

$$\phi^1 := \gamma(z_1, g_1 v), \quad (2.36)$$
$$\phi_0^1 := \gamma_0(z_1, g_1 v). \quad (2.37)$$

Adding up $\Sigma_{j>k} D_{3kj} D_{4kj}$ yields 1 due to $||\gamma(u,v)|| < 1$ [19, 22] and summing $\Sigma_{k\geq 0} D_{1k} D_{2k}$ constructs:

$$\Psi(u,v) = \phi_0^1 - \phi_0 + \frac{\xi \cdot \gamma_0}{1 - \xi\gamma}\left(\phi^1 - \phi\right). \tag{2.38}$$

Finally, we have

$$\Phi_M = \mathfrak{D}_{(u,v)}^{(M,M)}\left[\phi_0^1 - \phi_0 + \frac{\xi \cdot \gamma_0}{1 - \xi\gamma}\left(\phi^1 - \phi\right)\right]. \tag{2.39}$$

∎

The mean for the number of customers in the platform $A_\mu$ could be found from (2.39) and (2.41)-(2.42):

$$\mathbb{E}[A_\mu] = \mathbb{E}[\mathbb{E}[A_\mu|B_{\nu-1}]] = \lim_{z_1 \to 1}\sum_{k\geq 0}\left\{\frac{d}{dz_1}\Phi_M(1,1,z_1,1,1) p_k^b\right\} \tag{2.43}$$

where

$$p_k^b = \left(\frac{1}{k!}\right)\frac{\partial^k}{\partial g_0^k}\Phi_M(1,1,1,g_0,1)\bigg|_{g_0 \to 0}, \tag{2.44}$$

and

$$\Phi_M(1,1,1,g_0,1) = \mathfrak{D}_{(u,v)}^{(M,M)}\left[\phi_0^1 - \phi_0 + \frac{\xi \cdot \gamma_0}{1 - \xi\gamma}\left(\phi^1 - \phi\right)\right], \tag{2.45}$$

where, from (2.32)-(2.37),

$$\gamma := \gamma(u, g_0 v), \tag{2.46}$$
$$\gamma_0 := \gamma_0(u, g_0 v), \tag{2.47}$$

$$\phi := \gamma(u, v), \tag{2.48}$$
$$\phi_0 := \gamma_0(u, v), \tag{2.49}$$

$$\phi^1 := \gamma(1, v), \tag{2.50}$$
$$\phi_0^1 := \gamma_0(1, v). \tag{2.51}$$

**2.2. Open Two-sided Platform Model**

This innovative two-sided platform model could have several variants. One of variants is the two-sided platform has an unlimited capacity. Practically speaking, it is not possible that a real system has unlimited capacity or resource of the platform. Unlike the original stochastic two-sided platform model, the capacity of an open two-sided platform is relatively a large number and both sides are not interrelated. Recalling from (2.1)-(2.14), the initial setup of the open model is the same as the previous section. The open two-sided platform means that the capacity of the platform is almost unlimited and the joint functional of the two-sided platform becomes:

$$\Phi_M = \lim_{M \to \infty^-} \Phi_M(\xi, z_0, z_1, g_0, g_1). \tag{2.52}$$

Since the customer flow in a two-sided platform is given by the supplier flow from (2.15), the joint functional is modified only for the customer side to represent the open two-sided platform as follows:

$$\varphi(b; z_0, z_1) = \mathbb{E}\left[z_0^{A_{\mu(b)-1}} \cdot z_1^{A_{\mu(b)}} \middle| B_{\nu-1} = b\right], \tag{2.53}$$

which is given by the flows from the supplier side. From (2.41)-(2.43), the PGFs (Probability Generating Functions) of $A_\mu$, $A_{\mu-1}$ are as follows:

$$\varphi(b; z_0, z_1) = \gamma_a^0(z_0 z_1)\gamma_a(z_0 z_1)^{\mu(b)-1}\gamma_a(z_1), \tag{2.54}$$

where

$$\gamma_a(z) = \delta(\lambda_a(1-z)), \gamma_a^0(z) = \mathbb{E}[z^{A_0}] = \delta_0(\lambda_a(1-z)) \tag{2.55}$$

and the *exit index* of the customer side is:

$$\mu(b) := \mathbb{E}[\mu | B_\nu = b] \simeq \left\lfloor \frac{b - b_0}{\lambda_b \widetilde{\delta}} \right\rfloor, \tag{2.56}$$

($B_0$ is assumed as the fixed value $b_0 = \mathbb{E}[B_0]$). It is noted that all index values $\nu - 1$ and $\nu$ are almost the same (i.e., $\nu - 1 \simeq \nu$) when the resource of the supplier is almost unlimited (i.e., $0 \ll M < \infty$ and $0 \ll B_{\nu-1} \simeq B_\nu < \infty$). The mean number of customers $\mathbb{E}[A_\mu]$ in the two-sided platform could be found as follows:

$$\mathbb{E}\left[A_{\mu(B_{\nu-1})}\right] = \mathbb{E}\left[\mathbb{E}\left[A_{\mu(B_{\nu-1})} | B_{\nu-1}\right]\right], \tag{2.57}$$

where, from (2.15) and (2.56),

$$\mathbb{E}\left[A_{\mu(b)} | B_\nu = b\right] = \frac{d}{dz_1}\varphi(b; 1, z_1)\bigg|_{z_1 \to 1} = \lambda_a \widetilde{\delta_0} + \mu(b) \cdot \lambda_a(b)\widetilde{\delta}, \tag{2.58}$$

From (2.56) and (2.58), we have

$$\mathbb{E}\left[A_{\mu(B_{\nu-1})}\right] = \sum_{k=0}^{M} \left\{\lambda_a \widetilde{\delta_0} + \lambda_a(k - b_0)\right\} p_k^b \tag{2.59}$$

and

$$p_k^b = P\{B_\nu = k\} = \left(\frac{1}{\boldsymbol{A}}\right)\left\{\frac{\left(\lambda_a \widetilde{\delta}\right)^k}{k!} e^{-\left(\lambda_a \widetilde{\delta}\right)}\right\}, \tag{2.60}$$

$$\boldsymbol{A} = \sum_{k=0}^{M} \frac{\left(\lambda_a \widetilde{\delta}\right)^k}{k!} e^{-\left(\lambda_a \widetilde{\delta}\right)}. \tag{2.61}$$

Although the open platform model is relatively simple and easy to be analytically solvable, the open case is very practical for some specific cases.

## 2.3. Memoryless Two-sided Platform

It is assumed that the observation process has the memoryless properties which might be a special condition but very practical for actual implementation for analyzing the two-sided platform. It implies that the flow of suppliers does not contain any past information of the supplier side. We can find explicit solutions of $p_k^b$ and $\mathbb{E}[A_\mu]$ to build a proper payoff function for the optimization. Recall from (2.20), the $\mathcal{D}$-operator for single variable is defined as follows:

$$H(u) = \mathcal{D}_x[w(x)](u) := (1-u)\sum_{x \geq 0} w(x)u^x, \tag{2.62}$$

$$\mathcal{D}_{(x,y)}[w_1(x)w_2(y)](u,v) := (1-u)(1-v)\sum_{x \geq 0}\sum_{y \geq 0} w_1(x)w_2(y)u^x v^y$$
$$= \mathcal{D}_x[w_1(x)]\mathcal{D}_y[w_2(y)],$$

then

$$w(x,y) = \mathfrak{D}_{(u,v)}^{(x,y)}[\mathcal{D}_{(x,y)}\{w(x,y)\}], \tag{2.63}$$

$$w_1(x)w_2(y) = \mathfrak{D}_u^x[\mathcal{D}_x\{w_1(x)\}]\mathfrak{D}_v^y[\mathcal{D}_y\{w_2(y)\}], \tag{2.64}$$

where $\{f(x), (f_1(x)f_2(y))\}$ are a sequence, with the inverse (3.1) and

$$\mathfrak{D}_u^m(\bullet) = \begin{cases} \frac{1}{m!} \lim_{u \to 0} \frac{\partial^m}{\partial u^m} \frac{1}{(1-u)}(\bullet), & m \geq 0, \\ 0, & \text{otherwise,} \end{cases} \tag{2.65}$$

and

$$\mathfrak{D}_{(u,v)}^{(m,n)}[H_1(u)H_2(v)] = \mathfrak{D}_u^m[H_1(u)]\mathfrak{D}_v^n[H_2(v)]. \tag{2.66}$$

The function $\mathbb{E}[g_0^{B_{\nu-1}}]$ is the PGF (Probability Generating Function) of the number of suppliers at the prior moment of exceeding more than the platform capacity $M$. The probability distribution of $B_{\nu-1}$ could be found from (2.44) after obtaining the PGF from (2.42):

$$\mathbb{E}[g_0^{B_{\nu-1}}] = \Phi_M(1,1,1,g_0,1) = L_1 + \sum_{i \geq 0}\{Q_i^a - Q_{i+1}^a\}Q_i^b, \tag{2.95}$$

where

$$L_1 = \frac{\beta_b^0\left(1-(\alpha_b^0)^{M+1}\right)}{1-\alpha_b^0} - \frac{\beta_a^0 \beta_b^0 \left(1-(\alpha_a^0)^{M+1}\right)\left(1-(\alpha_b^0)^{M+1}\right)}{(1-\alpha_a^0)(1-\alpha_b^0)}, \tag{2.96}$$

$$Q_i^a = \beta_a^0(\beta_a)^i \sum_{j_a \geq 0} \Xi_i^M\left(j_a, \frac{\alpha_a}{\alpha_a^0}, \alpha_a^0\right), \tag{2.97}$$

$$Q_i^b = \left(\frac{\beta_a^0(\beta_a)^{i+1}}{\alpha_b - \alpha_b^0 g_0}\right)\sum_{j_2 \geq 0}\left\{(\alpha_b)\Xi_i^M(j_2, g_0, \alpha_b) - (\alpha_b^0 g_0)\Xi_i^M\left(j_2, \frac{\alpha_b}{\alpha_b^0}, \alpha_b^0 g_0\right)\right\}, \tag{2.98}$$

$$\Xi_i^m(j,x,y) = \binom{i-1+j}{j} x^j \left(\frac{y^j - y^{m+1}}{1-y}\right), |x| \leq 1, |y| \leq 1, \tag{2.99}$$

## III. STOCHASTIC OPTIMIZATIONS FOR TWO-SIDED PLATFORM

### 3.1. Two-sided Platform Payoff Function Design

The input ratio between two sides is a vital factor for a cost related function which includes the cost, the revenue and the payoff functions. The payoff function (i.e., a revenue function or a profit function) could be constructed based on the attraction factor $\alpha$ that indicates the customer coverage (i.e., the submissive side) by the supplier side (i.e., the dominant side). Practically, a typical range of the attraction factor $\alpha$ is up to 10 (i.e., $\alpha \in [1,10]$) although it could reach up to the infinity ($\alpha<\infty$). The optimal value of $\alpha$ towards to maximize the payoff function $\mathfrak{S}$ of the two-sided platform and the payoff function is:

$$\mathfrak{S}(\alpha; M)$$

$$= \mathbb{E}\left[c_0 \cdot A_\mu - c_1 \cdot \left\{B_\mu \cdot \mathbf{1}_{\{A_\mu+B_{\mu-1}\leq M\}} + (M - A_\mu) \cdot \mathbf{1}_{\{A_\mu+B_{\mu-1}>M\}}\right\}\right]$$

$$= c_0 \mathbb{E}[A_\mu] - c_1 \mathbb{E}\left[B_\mu \cdot \mathbf{1}_{\{A_\mu+B_{\mu-1}\leq M\}}\right] \quad (3.1)$$

$$- c_1\left\{M \cdot \mathbb{E}\left[\mathbf{1}_{\{A_\mu+B_{\mu-1}>M\}}\right] - \mathbb{E}\left[A_\mu \cdot \mathbf{1}_{\{A_\mu+B_{\mu-1}>M\}}\right]\right\},$$

and the best value $\alpha^*$ could be found as follows:

$$\alpha^* = \left\{\exists \alpha \in [1,\infty) : \frac{\partial \mathfrak{S}(\alpha; M)}{\partial \alpha} = 0\right\}, \quad (3.2)$$

where $c_0$ is the unit sales price from the customer side and $c_1$ is the unit cost price based on the suppliers. If $c_0 \geq 0$ and $c_1 \leq 0$, the platform earns the money from both sides. If $c_0 \leq 0$ and $c_1 \geq 0$, the platform spend the money from both sides. Although a payoff function could be constructed based on the attract factor in this section, other variables besides of the attraction factor could also construct payoff functions in different perspectives.

### 3.2. Two-sided Platform Optimization Setup And Visualization

The payoff function from (3.1) could be computationally implemented after all required parameters which are explained on Table 1 and these parameters are determined by actual data gathering. The practical real-world case (i.e., nightclub) shall be introduced on the next section. Although the payoff function in this section is constructed based on the attraction factor and the platform capacity, the payoff function could be designed in a different way.

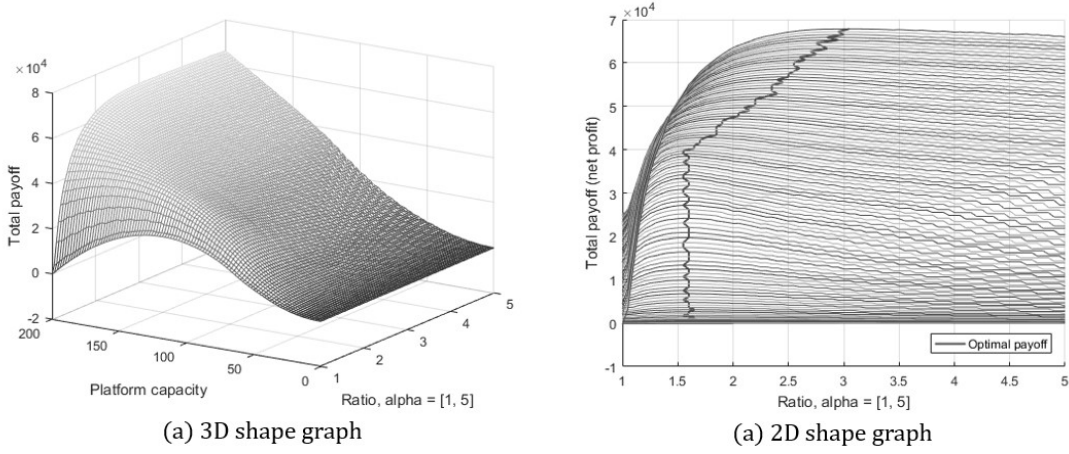

(a) 3D shape graph  (a) 2D shape graph

**Figure 1.** Visualizations of two-sided platform optimizations

The visualization of the payoff function for a two-sided platform model is shown in Fig. 1. The payoff (i.e., reward) function from (3.1) is a two-variable function which could be maximized based on two variables which are the attraction factor (i.e., the ratio between a supplier and a customer sides) and the platform capacity. It is noted that a two-variable function such as a payoff function of a two-sided platform could be visualized either as a 3D shape (a) or a 2D shape (b) in Fig.1. A proper shape depends on the visiualization objectives of a payoff function.

### 3.3. Properties of Open Two-sided Platforms

As it has been mentioned on Section 2, an open two-sided platform is the two-sided platform (almost) without the limitation of the capacity. From (3.1) the payoff function is revised as follows:

$$\lim_{M \to \infty} \mathfrak{S}(\alpha; M) = c_0 \cdot \mathbb{E}[A_\mu] - c_1 \cdot \mathbb{E}[B_\mu] = (c_0 - c_1 \cdot \alpha)\mathbb{E}[A_\mu], \quad (3.3)$$

where

$$\mathbb{E}[A_\mu] = \alpha \mathbb{E}[B_\mu]. \quad (3.4)$$

Due to unlimited capacity, more suppliers attract more customers and more profits. It indicates that there are no turning points of the payoff function but the attract factor $\alpha$ should be greater than a certain value to start making the profit of the two-sided platform. Hence, the optimized attraction factor of an open two-sided platform could be terminated as follows:

$$\alpha^* := \inf\left\{\alpha : \lim_{M \to \infty} \mathfrak{S}(\alpha; M) \geq 0\right\}, \quad (3.5)$$

From (3.4)-(3.5), we have

$$\alpha^* = \frac{c_1}{c_0}, \quad (3.6)$$

which indicates that the platform starts making benefits when the attraction factor is greater than $\alpha^*$. It is noted that the moment of making benefits depends on the cost ratio between the cost $c_1$ and the reward $c_0$ (i.e., $\alpha \geq \frac{c_1}{c_0}$).

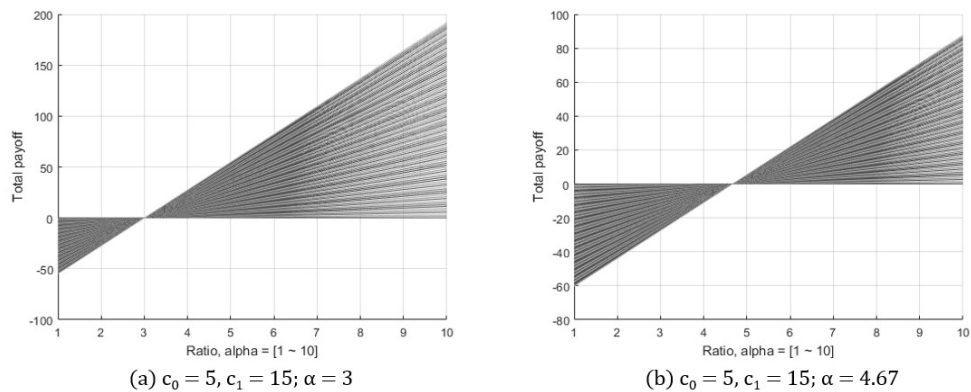

(a) $c_0 = 5, c_1 = 15; \alpha = 3$  (b) $c_0 = 5, c_1 = 15; \alpha = 4.67$

**Figure 2.** Open two-sided platform optimization examples